\documentclass[12pt]{article}
\usepackage{amssymb,amsmath,graphics}
\newtheorem{thm}[equation]{Theorem}
\newtheorem{lemma}[equation]{Lemma}

\newenvironment{pf}
{\begin{trivlist} \item \noindent{\sc Proof. }}
{{\hfill $\Box$} \end{trivlist}}
\newenvironment{pfof}[1]
{\begin{trivlist} \item \noindent{\sc Proof of\ #1. }}
{{\hfill $\Box$} \end{trivlist}}

\newcommand{\zd}{{\mathbb Z}^d}
\newcommand{\diam}{\mbox{\rm diam}}
\newcommand{\prob}{{\mathbf P}}
\newcommand{\expe}{{\mathbf E}}
\newcommand{\pss}{for $p$ sufficiently small}
\newcommand{\conn}[1]{\stackrel{#1}{\longleftrightarrow}}

\newcounter{mycount}
\newenvironment{mylist}{\begin{list}{{\rm (\roman{mycount})}}%
{\usecounter{mycount}\itemsep 0pt}}{\end{list}}

\newcommand{\dof}{\bf}

\begin{document}

\title{The Metastability Threshold for Modified Bootstrap Percolation
 in $d$ Dimensions}
\author{Alexander E. Holroyd}
\date{31 March 2006}

\maketitle
\renewcommand{\thefootnote}{}

\begin{abstract}
In the modified bootstrap percolation model, sites in the cube
$\{1,\ldots,L\}^d$ are initially declared active independently
with probability $p$.  At subsequent steps, an inactive site
becomes active if it has at least one active nearest neighbour in
each of the $d$ dimensions, while an active site remains active
forever.  We study the probability that the entire cube is
eventually active.  For all $d\geq 2$ we prove that as
$L\to\infty$ and $p\to 0$ simultaneously, this probability
converges to $1$ if $L=\exp^{d-1} \frac{\lambda+\epsilon}{p}$, and
converges to $0$ if $L=\exp^{d-1} \frac{\lambda-\epsilon}{p}$, for
any $\epsilon>0$.  Here $\exp^n$ denotes the $n$-th iterate of the
exponential function, and the threshold $\lambda$ equals $\pi^2/6$
for all $d$.
\footnote{\hspace{-2em} Funded in part by an NSERC (Canada)
Discovery Grant, and by MSRI (Berkeley USA)}
\footnote{\hspace{-2em}{\bf Address:} {\tt
holroyd(at)math.ubc.ca}. Department of Mathematics, University of
British Columbia, Vancouver, BC V6T 1Z2, CANADA}
\footnote{\hspace{-2em}{\bf Key words:} bootstrap percolation,
cellular automaton, metastability, finite-size scaling}
\footnote{\hspace{-2em}{\bf 2000 Mathematics Subject
Classifications:} Primary 60K35; Secondary 82B43}
\end{abstract}

\section{Introduction}

Let $\zd:=\{x=(x_1,\ldots,x_d):x\in {\mathbb Z}\}$ be the
$d$-dimensional integer lattice.  We call the elements of $\zd$
{\dof sites}.  Let
$e_1:=(1,0,\ldots,0),\ldots,e_d:=(0,\ldots,0,1)\in\zd$ be the
standard basic vectors.  For a set of sites $W\subseteq\zd$,
define
$$\beta(W):=W\cup\Big\{x\in\zd: \forall i=1,\ldots, d
\text{ we have } x+e_i\in W \text{ or } x-e_i\in W \Big\},$$ and
$$\langle W\rangle := \lim_{t\to\infty} \beta^t(W),$$
where $\beta^t$ denotes the $t$-th iterate of the function
$\beta$.  $\langle W\rangle$ is the final active set for the {\dof
modified bootstrap percolation model} starting with $W$ active.

Now fix $p\in(0,1)$ and let $X$ be a random subset of $\zd$ in
which each site is independently included with probability $p$.
More formally, denote by $\prob_p=\prob$ the product probability
measure with parameter $p$ on the product $\sigma$-algebra of
$\{0,1\}^{\zd}$, and define the random variable $X$ by
$X(\omega):=\{x\in\zd: \omega(x)=1\}$ for $\omega\in\{0,1\}^{\zd}$.
A site $x\in\zd$ is said to be {\dof occupied} if $x\in X$.

We say that a set $W\subseteq\zd$ is {\dof internally spanned}, or
{\dof i.s.}, if $\langle X\cap W\rangle=W$ (that is, if the model
restricted to $W$ fills $W$ up). For a positive integer $L$ we
define the $d$-dimensional {\dof cube} of side $L$ to be
$$Q^d(L):=\{1,\ldots ,L\}^d.$$
For convenience we also write $Q^d(L)=Q^d(\lfloor L\rfloor)$ when
$L$ is not an integer. We define
$$I^d(L)=I^d(L,p):=\prob_p\big(Q^d(L) \text{ is internally spanned}\big).$$
Let $\exp^n$ denote the $n$-th iterate of the exponential
function.

\begin{thm}
\label{main} Let $d\geq 2$ and $\epsilon>0$.  For the modified
bootstrap percolation model, as $p\to 0$ we have
\begin{mylist}
\item $\displaystyle
I^d\bigg(\exp^{d-1}\frac{\lambda+\epsilon}{p},\;p\bigg)\to 1;$ \item
$\displaystyle
I^d\bigg(\exp^{d-1}\frac{\lambda-\epsilon}{p},\;p\bigg)\to 0;$
\end{mylist}
where
$$\lambda=\frac{\pi^2}{6}.$$
\end{thm}

\subsubsection*{Remarks}

The case $d=2$ of Theorem \ref{main} was proved in \cite{h-boot}.
The modified bootstrap percolation model considered here is a
minor variant of the {\dof standard bootstrap percolation model},
which is defined in the same way except replacing the function
$\beta$ with
$$\beta'(W):=W\cup\Big\{x\in\zd: \#\{y\in W:\|y-x\|_1=1\}\geq
d\Big\},$$
(so a site becomes active if it has at least $d$ active
neighbours).
 In the case $d=2$, the analogue of Theorem \ref{main}
was proved for the standard bootstrap percolation model in
\cite{h-boot}; in this case the threshold $\lambda$ becomes
$\pi^2/18$.  Similar results were obtained for a further family of
two-dimensional models in \cite{h-liggett-romik}.  The present
work is the first proof of the existence of a sharp threshold
$\lambda$ for a bootstrap percolation model in 3 or more
dimensions; in addition we determine the value $\pi^2/6$. The
analogue of Theorem \ref{main} but with two {\em different}
constants $c_1,c_2$ in place of
$\lambda+\epsilon,\lambda-\epsilon$ was proved earlier in
\cite{aiz-leb} ($d=2$), \cite{cerf-c} ($d=3$) and \cite{cerf-m}
($d\geq 4$).  These works apply to the standard model (among
others), but can easily be adapted to the modified model
considered here.  It is a fascinating open problem to prove the
existence of a sharp threshold for the standard model in $3$ or
more dimensions.

In \cite{balogh-bollobas-sharp} a different kind of ``sharpness"
is proved, by a general method, for various models including
standard and modified bootstrap percolation: writing
$p_\alpha=p_\alpha(L)$ for the value such that
$I(L,p_\alpha)=\alpha$, then
$p_{1-\epsilon}-p_\epsilon=o(p_{1/2})$ as $L\to\infty$, with a
certain explicit bound. (However this result says nothing about
the behaviour of $p_{1/2}$ as a function of $L$).  Similar results
with the roles of $p$ and $L$ exchanged may be obtained using the
methods of \cite{aiz-leb}.

There have been numerous other beautiful rigorous contributions to
the study of bootstrap percolation models, initiated by
\cite{van-enter}.  For example see the references in
\cite{cerf-m},\cite{h-boot}.

Bootstrap percolation models have important applications, both
directly and as tools in the study of more complicated systems
(see for example the references in \cite{cerf-m},\cite{h-boot}).
The models have been extensively studied via simulation, and it is
a remarkable fact that the resulting asymptotic predictions often
differ greatly from rigorous asymptotic results, apparently
because the convergence as $(L,p)\to(\infty,0)$ is extremely slow.
See \cite{schonmann-boot},\cite{h-boot} for examples.  In the case
of the modified bootstrap percolation model in $d=2$, the value
$0.47\pm 0.02$ for $\lambda$ was predicted numerically in
\cite{a-s-a}, whereas the rigorous result from \cite{h-boot} is
$\lambda=\pi^2/6=1.644934\cdots$.  It would be worthwhile to
compare simulations with our rigorous result that
$\lambda=\pi^2/6$ for $d\geq 3$.  It is of interest to understand
this slow convergence phenomenon in more detail, and it is
relevant to applications: a typical physical system might have
$L^d\approx 10^{20}$ particles, which is much larger than current
computer simulations allow, but potentially not large enough to
exhibit a threshold close to the limiting value.  See \cite{glbd}
for an interesting partly non-rigorous investigation of some these
issues.

\subsubsection*{Proof outline}

The proof of Theorem \ref{main} is by induction on the dimension.
The base case $d=2$ is provided by the results in \cite{h-boot}.
(The proof in \cite{h-boot} is quite involved, and very specific
to the 2-dimensional model). It is interesting that the constant
$\lambda=\pi^2/6$ enters only here.  The inductive step follows
closely the pioneering work of \cite{cerf-c},\cite{cerf-m}, although since
our result is more precise we need to be more careful with the
estimates.  As in \cite{cerf-m}, the case $d=3$ is the most
delicate.

The proof of the lower bound in Theorem \ref{main}(i) is
relatively straightforward, and many of the ideas were already
present in \cite{schonmann-boot}.  The fundamental observation is
that if a cube is already entirely active, then the sites lying on
its faces evolve according to the modified bootstrap percolation
model in $d-1$ dimensions.  Hence, by the inductive hypothesis, a
cube of size $L$ is likely to be internally spanned if it contains
some internally spanned cube of size
$m=\exp^{d-2}\frac{\lambda+\epsilon}{p}$ (sometimes called a
``critical droplet" or ``nucleation centre"), because such a cube
will grow forever from its faces.  For $d\geq 3$, straightforward
arguments show that such a cube is internally spanned with
probability at least (roughly) $e^{-m}$, and so in order to
internally span the larger cube we should take approximately
$L>1/(e^{-m})=\exp^{d-1}\frac{\lambda+\epsilon}{p}$, completing
the induction.

The proof of the upper bound in Theorem \ref{main}(ii) is more
challenging, and is based on the more subtle construction originating
in \cite{cerf-c}.  The idea is to find an upper bound on the
probability that a cube of size
$m=\exp^{d-2}\frac{\lambda-\epsilon}{p}$ has a left-right crossing in
its final configuration.  (Such a crossing plays the role of a
nucleation centre in this bound).  The proof proceeds by dividing this
cube into ``slices", and running the $(d-1)$-dimensional model in
each, to produce a configuration which dominates the $d$-dimensional
model.  By the inductive hypothesis, the probability that a slice
becomes fully active is small, and, where a slice does not become
fully active, its final configuration resembles subcritical
percolation.  Hence a connection of length $m$ has probability at most
(roughly) $e^{-m}$, and again we can complete the induction.  A key
point in the present proof is that for the modified model, we can use
slices of thickness 1, whereas in \cite{cerf-c},\cite{cerf-m} (for the
standard model) it was necessary to use slices of thickness 2, and to
replace the parameter $p$ with $2p$.  Changing $p$ in this way makes
it impossible to obtain matching upper and lower bounds, and it is for
this reason that our method cannot be adapted directly to prove an
analogous result for the standard bootstrap percolation model.
Another difference in the proof here as compared with
\cite{cerf-c},\cite{cerf-m} is that (in the case $d=3$) we need to
carefully balance the probabilities of fully active slices with those
of percolation connections.  Equation (\ref{cases-bound}) is the heart
of this calculation.

\subsubsection*{Notation and conventions}

It will be convenient to consider lower-dimensional versions of
the model running on subsets $\zd$.  Let
$\delta\in\{1,\ldots,d\}$. We define the $\delta$-dimensional cube
$$Q^\delta(L):=\{1,\ldots ,L\}^\delta\times\{1\}^{d-\delta}\subseteq\zd.$$
By a {\dof copy} of $Q^\delta(L)$ we mean an image of
$Q^\delta(L)$ under any isometry of $\zd$.  For a set
$W\subseteq\zd$, which will always be a subset of some copy of a
$\delta$-dimensional cube, we define
$$\beta_\delta(W):=W\cup\Big\{ x\in\zd:
\#\big\{i:\{x+e_i,x-e_i\}\cap W \neq\emptyset\big\}\geq
\delta\Big\},$$ and $\langle W\rangle_\delta := \lim_{t\to\infty}
\beta_\delta^t(W).$ We say that $W$ is $\delta${\dof -internally
spanned} if $\langle X\cap W\rangle_\delta=W$.  Let
$I^\delta(L):=\prob_p(Q^\delta(L) \text{ is $\delta$-i.s.})$, and
note that this is consistent with the earlier definition.

Theorem \ref{main} involves an asymptotic statement as $p\to 0$
with $d$ and $\epsilon$ fixed.  Many of the inequalities used in
the proof will be valid ``\pss", be which we mean for all $p$ less
than some $q=q(d,\epsilon)>0$, whose value may vary from one
instance to another. In some of following proofs we use
$C_1,C_2,\ldots$ to denote constants in $(0,\infty)$ which may
depend on $d$ and $\epsilon$, but not on $p$.

\section{Lower Bound}

\begin{lemma}
\label{seeds} For any $d\geq 3$ and for $a,\epsilon>0$, if $p$ is
sufficiently small (depending on $d,a,\epsilon$) then
$$I^d\bigg(\exp^{d-2}\frac{a}{p}\bigg)\geq 1/\exp^{d-1}\frac{a+\epsilon}{p}.$$
\end{lemma}

\begin{figure}
\centering \resizebox{!}{1in}{\includegraphics{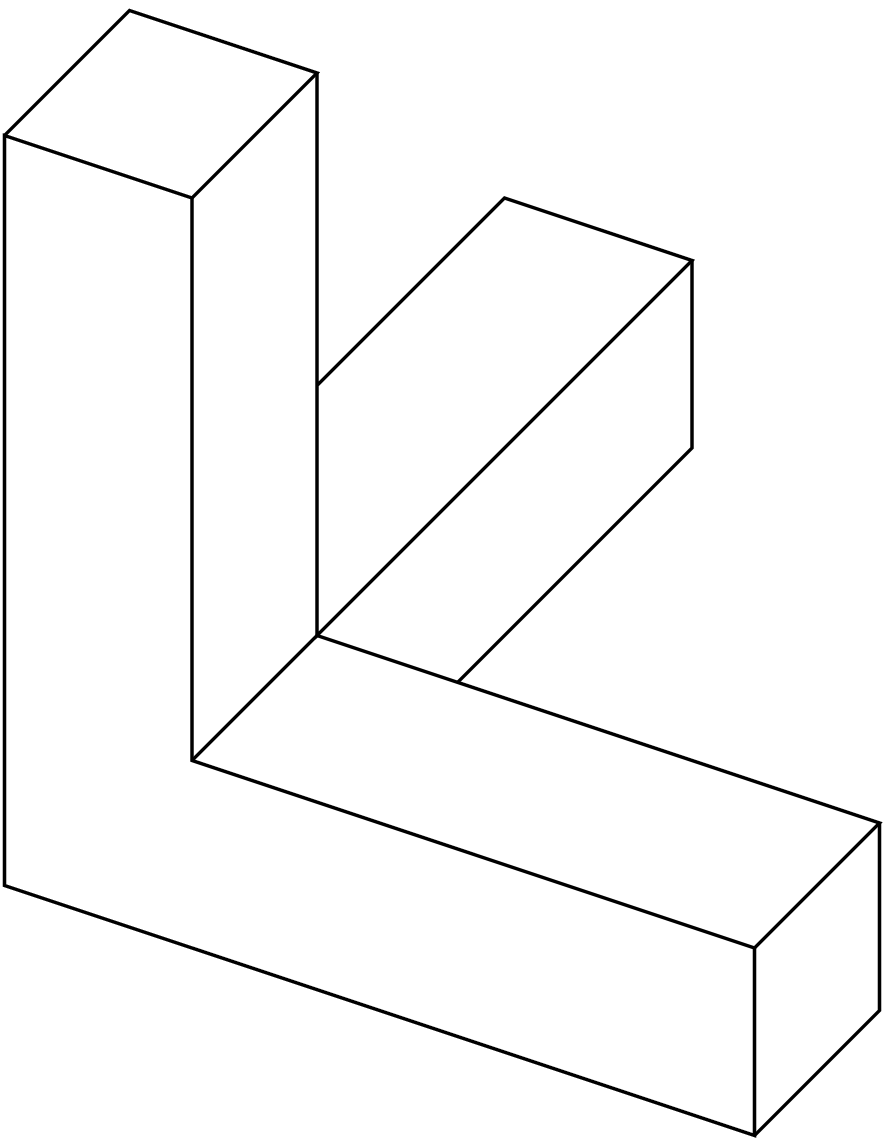}}
\caption{The occupied set in event $A$.}
\label{ell}
\end{figure}

\begin{pf}
If $d\geq 4$, note that a cube is internally spanned if {\em all}
of its sites are occupied.  Therefore
$$I^d\bigg(\exp^{d-2}\frac{a}{p}\bigg)\geq
p^{\big[\exp^{d-2}(a/p)\big]^d}\geq 1/\exp^{d-1}\frac{a+\epsilon}{p}$$
\pss. (To check the second inequality, take three successive
logarithms of the reciprocal of both sides).

The case $d=3$ is a little more delicate.  Write $L=\lfloor
e^{a/p}\rfloor$ and $k=\lfloor p^{-3}\rfloor$ (so $k\ll L$ for $p$
sufficiently small). Let $A$ be the event that every site having
two of its coordinates in $\{1,\ldots,k\}$ and one coordinate in
$\{1,\ldots,L\}$ is occupied (see Figure \ref{ell}).
Let $B$ be the event that every
copy of $Q^1(k)$ in $Q^3(L)$ contains at least one occupied site.
It is straightforward to check that if $A$ and $B$ both occur then
$Q^3(L)$ is internally spanned.  Since $A$ and $B$ are increasing
events, the Harris-FKG inequality (see e.g.\ \cite{g2}) yields
$$I^3(L)\geq \prob(A)\prob(B).$$
We now estimate
$$\prob(A)\geq p^{3Lk^2} \geq 1/\exp^2 \frac{a+\epsilon}{p}$$
\pss\ (to check the second inequality, take two logarithms of the
reciprocals), while
$$\prob(B)\geq 1-3L^3(1-p)^k\geq 1-3\exp(\lfloor a/p\rfloor -\lfloor p^{-2}
\rfloor)\geq e^{-1}$$
\pss.
Therefore we have
$$I^3(L)\geq 1/\exp^2 \frac{a+2\epsilon}{p}$$
\pss, and since $\epsilon$ was arbitrary this proves the result for $d=3$.
\end{pf}

The following result from \cite{aiz-leb} states that $I^d(L)$ increases rapidly with $L$ once it is large enough.
\begin{lemma}
\label{perc} For each $d\geq 1$ there exist $c=c(d)<1$ and
$C=C(d)<\infty$ such that provided $I^d(\ell)\geq c$, we have for
all $L\geq \ell$ that
$$I^d(L)\geq 1-Ce^{-L/\ell}.$$
\end{lemma}

\begin{pf}
See \cite{aiz-leb},\cite{schonmann-boot} or \cite{cerf-m}.  The idea is
to divide $Q^d(L)$ into disjoint or nearly-disjoint copies of
$Q^d(\ell)$. If $Q^d(L)$ is not i.s.\ then it is crossed by a path
of non-i.s.\ copies of $Q^d(\ell)$.  The probability of this event
can be estimated using standard percolation methods.
\end{pf}


\begin{pfof}{Theorem \ref{main}(i)}
The proof is by induction on $d$.  The statement of the theorem
holds in the case $d=2$ by Theorems 4 and 1(i) of \cite{h-boot}.
Now let $d\geq 3$, and suppose that for all $\delta=2,\ldots ,d-1$
we have
\begin{equation}
\label{ind-hyp} \text{for every $\epsilon>0$,}\quad
I^\delta\Big(\exp^{\delta-1}\frac{\lambda+\epsilon}{p}\Big)\to 1
\text{ as }p\to 0.
\end{equation}
We shall deduce that (\ref{ind-hyp}) holds for $\delta=d$ also.

Fix $\epsilon>0$.
We first claim that
\begin{equation}
\label{big-seeds} I^d\big(\exp^{d-2}p^{-2}\big)\geq
1/\exp^{d-1}\frac{\lambda+4 \epsilon}{p}
\end{equation}
\pss .
To prove this, write
$$\ell=\Big\lfloor \exp^{d-2}\frac{\lambda+\epsilon}{p}\Big\rfloor;
\quad a=\Big\lfloor
\exp^{d-2}\frac{\lambda+2\epsilon}{p}\Big\rfloor; \quad
b=\Big\lfloor \exp^{d-2}p^{-2}\Big\rfloor.
$$
By Lemma \ref{seeds} we have
\begin{equation}
\label{small-seed}
I^d(a)\geq 1/\exp^{d-1}\frac{\lambda+3\epsilon}{p}.
\end{equation}
We will deduce the claimed lower bound on $I^d(b)$ using the fact
that an internally spanned cube will grow if each of its faces (of
all possible dimensions) is internally spanned in the model of the
appropriate lower dimension.

More precisely, for $L\geq 1$ and a proper subset $S\subsetneq
\{1,\ldots,d\}$, define the {\dof face}
$$F_S(L):=\big\{x\in\zd: x_i\in[1,L] \;\forall\; i\in S,
 \text{ and } x_i=L+1 \;\forall\; i\notin S\big\}.$$
Thus $F_S(L)$ is a copy of $Q^{|S|}(L)$, and we have the disjoint union
$$Q^d(L+1)=Q^d(L)\sqcup \bigsqcup_{S\subsetneq \{1,\ldots,d\}} F_S(L).$$
It is straightforward to check that if $Q^d(L)$ is $d$-i.s.\ and
the face $F_S(L)$ is $|S|$-i.s.\ for every $S\subsetneq
\{1,\ldots,d\}$ then $Q^d(L+1)$ is $d$-i.s.\ Hence we have
\begin{equation}
\label{growth}
I^d(b)\geq I^d(a) \;\prob(G_a^b),
\end{equation}
where $G_a^b$ is the event that $F_S(j)$ is $|S|$-i.s.\ for every
$j\in[a,b)$ and every $S\subsetneq \{1,\ldots,d\}$.

In order to bound $\prob(G_a^b)$, first note that we may take $p$
sufficiently small that
\begin{equation}
\label{all-c}
I^\delta(\ell)\geq c(\delta) \quad\text{ for all }\delta\in [1,d-1],
\end{equation}
where $c(\delta)$ is as in Lemma \ref{perc}.  (The case
$\delta=d-1$ follows directly from (\ref{ind-hyp});  the cases
$\delta\in[2,d-2]$ follow from (\ref{ind-hyp}) by an additional
application of Lemma \ref{perc};  the case $\delta=1$ is trivial.)
Therefore writing $C'=\max_{\delta\in[1,d-1]} C(\delta)$, Lemma
\ref{perc} yields
$$
\prob(G_a^b) \geq  1-2^d C'\sum_{j=a}^{b-1} e^{-j/\ell}
\geq  1-2^{d+1} C' \ell e^{-a/\ell}
\geq  e^{-1}
$$
\pss.
Combining this with (\ref{small-seed}),(\ref{growth}) proves
the claim (\ref{big-seeds}).

Now write $$L=\Big\lfloor \exp^{d-1}
\frac{\lambda+5\epsilon}{p}\Big\rfloor.$$ Let $E$ be the event
that $Q^d(L)$ contains some $d$-i.s.\ copy of $Q^d(b)$, and let
$F$ be the event that for each $\delta\in[1,d-1]$, every copy of
$Q^\delta(b)$ in $Q^d(L)$ is $\delta$-i.s.\  It is straightforward
to check that if $E$ and $F$ both occur then $Q^d(L)$ is $d$-i.s.
Hence by the Harris-FKG inequality,
\begin{equation}
\label{fkg}
I^d(L)\geq \prob(E)\;\prob(F).
\end{equation}

By tiling $Q^d(L)$ with disjoint copies of $Q^d(b)$, we have
\begin{eqnarray*}
\prob(E)&\geq& 1-\big(1-I^d(b)\big)^{\lfloor L/b\rfloor^d}\\
&\geq & 1-\exp\big[-I^d(b)(L/b)^d\big]
\end{eqnarray*}
But by (\ref{big-seeds}), \pss\ we have
$$I^d(b)(L/b)^d\geq I^d(b) L/b\geq
\frac {\Big\lfloor\exp^{d-1}
\frac{\lambda+5\epsilon}{p}\Big\rfloor}
{\exp^{d-1}\frac{\lambda+4\epsilon}{p}  \exp^{d-2}p^{-2}} \geq
\exp^{d-1} \frac{\lambda+4.9\epsilon}{p} \to \infty$$ as $p\to 0$,
hence $\prob(E)\to 1$.

On the other hand, again taking \pss\ to satisfy (\ref{all-c}), we
have by Lemma \ref{perc},
$$\prob(F)\geq 1-2^d L^d C' e^{-b/\ell}.$$
But we have
$$\log(L^d e^{-b/\ell})\leq d\exp^{d-2}
\frac{\lambda+5\epsilon}{p}-\frac{\lfloor \exp^{d-2}p^{-2}\rfloor}
{\exp^{d-2}\frac{\lambda+\epsilon}{p}} \leq -\exp^{d-2} (p^{-2}/2)$$
\pss.  Hence $\prob(F)\to 1$ as $p\to 0$.

Thus by (\ref{fkg}) we have proved that
$I^d(\exp^{d-1}\frac{\lambda+5\epsilon}{p})\to 1$ as $p \to 1$.
Since $\epsilon$ was arbitrary this is (\ref{ind-hyp}) with
$\delta=d$, and the induction is complete.
\end{pfof}

\section{Upper Bound}

The main step in the proof of Theorem \ref{main}(ii) will Theorem
\ref{sub-crit} below, which states that within a cube of
appropriate size, the final configuration of the model resembles
(highly) subcritical percolation.

We call a set of sites $W\subseteq\zd$ {\dof connected} if it
induces a connected graph in the nearest-neighbour hypercubic
lattice.  A {\dof component} is a maximal connected subset.  For
sites $x,y\in\zd$ and a (random) set $W\subseteq\zd$ we write
``$x\conn{W}y$" for the event that $W$ has a component containing
$x$ and $y$.  (Note that $x\conn{W}x$ is equivalent to $x\in W$).
For sites $x,y\in Q^d(m)$, we define
$$f_m^d(x,y)=f_m^d(x,y,p):=
\prob_p\bigg(x \conn{\langle X\cap Q^d(m)\rangle}y\bigg).$$

\begin{thm}
\label{sub-crit} Let $d\geq 3$ and $\epsilon>0$.  Let
$$m=\bigg\lfloor \exp^{d-2} \frac{\lambda-\epsilon}{p} \bigg\rfloor,$$
where $\lambda=\pi^2/6$.
 There exist $\gamma=\gamma(d,\epsilon)>0$ and
$q=q(d,\epsilon)>0$ such that for all $p<q$ and all $x,y\in
Q^d(m)$,
$$f_m^d(x,y)\leq p^{\gamma(\|x-y\|_\infty+1)}.$$
\end{thm}
The ``+1" term in the exponent is important, since for the
induction we need a bound which is $o(1)$ as $p\to 0$ even in the
case $x=y$.

The following result from \cite{aiz-leb} is very useful.  The
{\dof diameter} of a set $W\subseteq\zd$ is $\diam\,
W:=\sup_{x,y\in W} \|x-y\|_\infty.$
\begin{lemma}
\label{aiz-leb} If $S\subseteq\zd$ is connected and internally
spanned then for every real $a\in[1, \diam\, S]$ there exists a
connected, internally spanned set $T\subseteq S$ with $\diam\,
T\in[a,2a]$.
\end{lemma}
\begin{pf}
See \cite{aiz-leb},\cite{h-boot} or \cite{cerf-m}.  The idea is to
realize the bootstrap percolation model by an iterative algorithm.
We keep track of a collection of disjoint i.s.\ sets $S_j$.  At
each step, if there is a site in $\langle \cup_j
S_j\rangle\setminus\cup_j S_j$ then we unite it with at most $d$
of the sets to form a new set. In the case of the modified model,
$\max_j(\diam\, S_j)$ is at most doubled at each step.
\end{pf}

\begin{pfof}{Theorem \ref{main}(ii)}
The case $d=2$ was proved in \cite{h-boot}.  Therefore fix $d\geq
3$ and $\epsilon>0$, and let
$$L=\Big\lfloor\exp^{d-1}{\frac{\lambda-\epsilon}{p}}\Big\rfloor
\quad\text{and}\quad
m=\Big\lfloor\exp^{d-2}{\frac{\lambda-\epsilon}{p}}\Big\rfloor.
$$
By Lemma \ref{aiz-leb}, if $Q^d(L)$ is internally spanned then it
contains some connected internally spanned set $T$ with $\diam\,
T\in[m/2,m]$.  This implies that there exist $v\in Q^d(L)$ and
$x,y\in v+Q^d(m)\subset Q^d(L)$ with $\|x-y\|\geq m/2$ such that
$x\conn{\big\langle X\cap (v+Q^d(m))\big\rangle} y$.  Hence by
Theorem \ref{sub-crit} we obtain
\begin{align}
\label{sub-set}
 I^d(L)\leq L^d m^d m^d p^{\gamma (m/2+1)}  \leq 1/\exp^{d-1}
\frac{\lambda-\epsilon}{p} \to 0 \quad\text{ as }p\to 0
\end{align}
as required. Here $L^d m^d m^d$ is a bound on the number of
choices for $v,x,y$, and the second inequality holds \pss, by
taking the logarithm thus:
\begin{align*}
\log\big(L^d m^{2d} p^{\gamma (m/2+1)}\big)&\leq
 \bigg(d-\frac{\gamma}{2}\log\frac{1}{p}\bigg)
 \exp^{d-2}\frac{\lambda-\epsilon}{p}
 +2d \exp^{d-3}\frac{\lambda-\epsilon}{p} \\
&\leq (-1)\exp^{d-2}\frac{\lambda-\epsilon}{p}.
\end{align*}
\end{pfof}

The proof of Theorem \ref{sub-crit} is by induction on the
dimension.  The key estimate is Lemma \ref{constr} below, for
which we need to define two more quantities.
Let
$$\chi^d_n=\chi^d_n(p):=\sum_{y\in Q^{d}(2n+1)}
 f_{2n+1}^{d}(z,y)$$
where $z:=(n+1,\ldots,n+1)$ is the site at the centre of
 $Q^{d}(2n+1)$.  (Thus $\chi^d_n$ is the expected volume of the component
 at $z$ in the final configuration of the model on $Q^{d}(2n+1)$).
For $n\leq m$ define
$$F_{m,n}^d=F_{m,n}^d(p):=
\prob_p\Big(\langle X\cap Q^d(m)\rangle \text{ has a component of
diameter $\geq n$}\Big).$$

\begin{lemma}
\label{constr} For any $n\leq m$ and $x,y\in Q^d(m)$ with
$\ell=\|x-y\|_\infty$ we have
\begin{multline*}
f_m^d(x,y)\leq  \sum_{k=0}^\ell \;\;
\sum_{0<i_1<\cdots<i_k<i_{k+1}=\ell+1} \\
\left(H(i_1\!+\!1)+\sum_{a=0}^m F_{m,n}^{d-1}\Big[1\!\wedge\! m^{d-1}
 H(i_1\!+\!a)\Big]\right)
 \prod_{j=1}^k F_{m,n}^{d-1}\Big[1\wedge m^{d-1}
 H(i_{k+1}-i_k)\Big]
\end{multline*}
where
$$H(r):=\frac 12 \sum_{s=r-1}^\infty \left(2 \chi^{d-1}_n \right)^s.$$
\end{lemma}
(Perhaps the easiest way to understand Lemma \ref{constr} is to read the proof as far as (\ref{case-1}),(\ref{case-2}), and look at Figures \ref{2poss},\ref{vw}.)

\begin{pf}
The following construction is based on that of
\cite{cerf-c},\cite{cerf-m}. Without loss of generality suppose
that $\|x-y\|_\infty=(y)_d-(x)_d$; if not we reorder the
coordinates and/or reverse the direction of the $d$th coordinate.
Write $(x)_d=u$ and $(y)_d=u+\ell$. Divide the cube $Q^d(m)$ into
the {\dof slices}
$$T_j:=Q^{d-1}(m)\times\{j+u\},\quad j=-u+1,\ldots ,m-u,$$
so that $x\in T_0$ and $y\in T_\ell$.
Let $Y_j:=\langle X \cap T_j\rangle_{d-1}$ be the final
configuration of the $(d-1)$-dimensional model restricted to
$T_j$.  Let
$$Z_j:=
\left\{%
\begin{array}{ll}
    T_j & \quad\text{if $Y_j$ has a component of diameter $\geq n$;} \\
    Y_j & \quad\text{otherwise.} \\
\end{array}%
\right.
$$
In the former case we say that the slice $T_j$ is {\dof full}.
Now let
$$Z:=\bigcup_{j=-u+1}^{m-u} Z_j.$$

The point of this construction is that $Z\supseteq \langle X\cap
Q^d(m)\rangle_d$.  To see this note that $Z_j\supseteq
Y_j\supseteq\langle X\cap Q^d(m)\rangle_d\cap T_j$; the latter
inclusion holds because running the $(d-1)$-dimensional model on
$T_j$ is equivalent to running the $d$-dimensional model with the
boundary condition that every site in $Q^d(m)\setminus T_j$ is
occupied -- hence it must result in a larger configuration in
$T_j$ than running the $d$-dimensional model on $X\cap Q^d(m)$.
(Note that the argument would not work in this form for the
standard bootstrap percolation model, since the boundary condition
adds two extra neighbours to each site in the slice.) Therefore
$$f_m^d(x,y)\leq \prob \Big(x\conn{Z} y\Big).$$

We shall bound the above probability by splitting the event up
according to which slices are full.  Let $I_1<\cdots<I_K$ be the
random indices of those slices among $T_1,\ldots,T_{\ell}$ that
are full. Also let ${\cal W}$ be the event than every path in $Z$
from $x$ to $y$ intersects some full slice among
$T_{-u},\ldots,T_0$, and let $-A$ be the index of the last full
slice among $T_{-u},\ldots,T_0$ (or $A=\infty$ if there is none).
Then
\begin{align}
\nonumber f_m^d(x,y)\leq & \sum_{k=0}^\ell \;\;
\sum_{0<i_1<\cdots<i_k<\ell+1} \;\;
\sum_{a=0,\ldots,m,\infty} \\
\label{case-1} &\;\;\bigg[ \prob \Big(x\conn{Z} y,\;
(I_1,\ldots ,I_K)=(i_1,\ldots ,i_k),\;A=a,\; {\cal W}^C\Big)\\
\label{case-2} &\;\;\;+\prob \Big(x\conn{Z} y,\; (I_1,\ldots
,I_K)=(i_1,\ldots ,i_k),\; A=a,\; {\cal W} \Big) \bigg].
\end{align}
(See Figure \ref{2poss} for an illustration).
\begin{figure}
\centering \resizebox{!}{1.5in}{\includegraphics{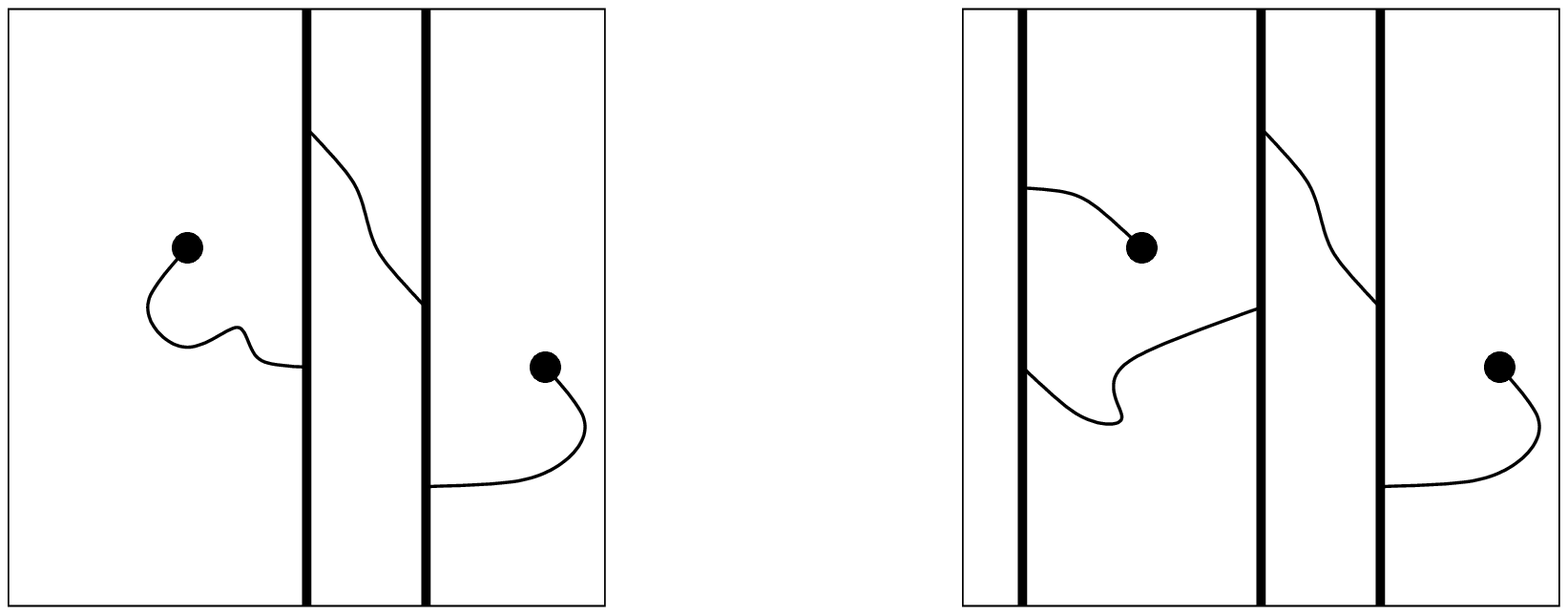}}
\caption{Two possibilities for the connection from $x$ to $y$
(black dots).  Vertical bars indicate full slices.  In the left
picture ${\cal W}$ does not occur; in the right picture it does.}
\label{2poss}
\end{figure}

Using independence of the slices, the probability (\ref{case-2}) above is
at most
\begin{multline}
\label{bound-1} \prob \big(T_{-a} \text{ is full}\big)\prob
\big({\cal E}(-a+1,i_1-1)\big)
\\
\times \prod_{j=1}^k \prob \big( T_{i_j} \text{ is full}\big)\prob
\big({\cal E}(i_j+1,i_{j+1}-1)\big),
\end{multline}
where we have written for convenience $i_{k+1}:=\ell+1$,
and where
$${\cal E}(i,i'):=\bigg\{T_i,\dots,T_{i'}\text{ are not full, and }
v\conn{\bigcup_{j=i}^{i'} Z_j} v'
\text{ for some }v\in T_i,\; v'\in T_{i'}\bigg\}
$$
(and taking ${\cal E}(i,i')$ to be an event of probability 1 if $i>i'$).

Similarly, the probability (\ref{case-1}) is
at most
\begin{multline}
\label{bound-2}
\prob \big({\cal G}(x,-a+1,i_1-1)\cap\{A=a\}\big)\\
\times \prod_{j=1}^k \prob \big( T_{i_j} \text{ is full}\big)\prob
\big({\cal E}(i_j+1,i_{j+1}-1)\big),
\end{multline}
where
$${\cal G}(x,i,i'):=\bigg\{T_i,\dots,T_{i'}\text{ are not full, and }
x\conn{\bigcup_{j=i}^{i'} Z_j} v' \text{ for some } v'\in T_{i'}\bigg\}.
$$

Next we bound the factors in (\ref{bound-1}),(\ref{bound-2}).
For any slice $T_j$ we have
\begin{equation}
\label{prob-full} \prob (T_j \text{ is full})=F^{d-1}_{m,n},
\end{equation}
so it remains only to bound the probabilities of ${\cal E}(i,i')$
and ${\cal G}(x,i,i')$.

Suppose that the event ${\cal E}(i,i')$ occurs.  Then there is a
nearest-neighbour path in $\bigcup_{j=i}^{i'} Z_j$ from some site
in $T_i$ to some site in $T_{i'}$.  Let $\alpha$ be such a path of
minimum length.  Define sites $v_1,w_1,v_2,w_2,\ldots,w_s$ along
the path as follows (see Figure \ref{vw} for an illustration). Let
$v_1\in T_i$ be the first site of $\alpha$. Given $v_1,\ldots
,v_t$, let $v_{t+1}$ be the first site after $v_t$ at which the
path enters a slice different from that of $v_t$.  Iterate this
until when we reach a site $v_s\in T_{i'}$. Let $w_t$ be the site
preceding $v_{t+1}$ in $\alpha$ for each $t< s$, and let $w_s\in
T_{i'}$ be the last site of $\alpha$.  Thus $\alpha$ consists of a
sequence of $s$ sub-paths $(v_1,\ldots,w_1),\;\ldots,\;
(v_s,\ldots,w_s)$, each one lying entirely within one slice, and
with $w_t,v_{t+1}$ being adjacent but in different slices.  (Note
however that two non-adjacent sub-paths may lie in the same
slice).
More precisely we have the following facts.  For each $t=1,\ldots ,s$:
\begin{mylist}
 \item $v_t,w_t\in T_{j(t)}$ for some $j(t)\in[i,i']$, with
$j(1)=i$ and $j(s)=i'$;
 \item $|j(t)-j(t+1)|=1$ and
$\|w_t-v_{t+1}\|_1=1$ for $t<s$;
 \item $\|v_t-w_t\|_\infty\leq n$;
 \item the component of $Z_{j(t)}$ at $v_t$ has diameter $\leq n$;
 \item $v_t\conn{Z_{j(t)}}w_t$ occurs;
 \item $v_t\conn{Z_{j(t)}}v_{t'}$ does not occur for any $t'\neq t$.
\end{mylist}
Properties (iii),(iv) hold because the slices $T_i,\ldots ,T_{i'}$
are not full, and property (vi) holds because we chose a path of
minimum length.
\begin{figure}
\centering \resizebox{!}{2in}{\includegraphics{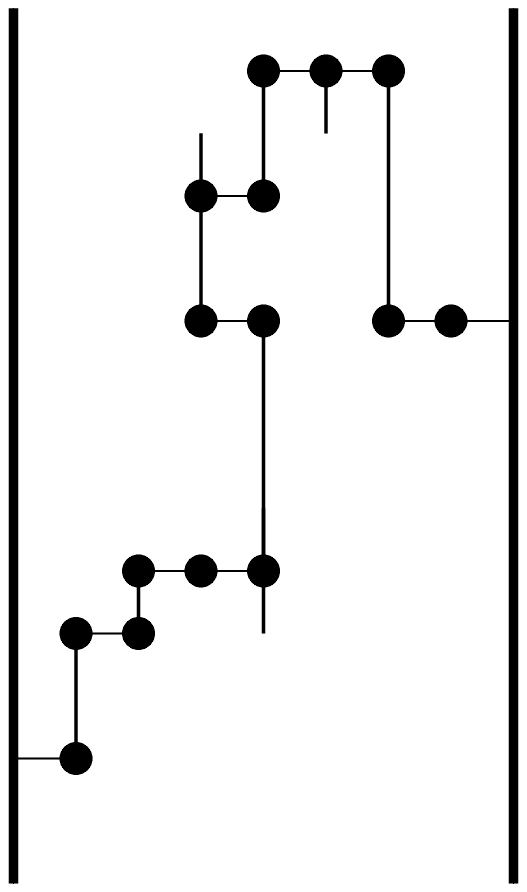}} \caption{An
illustration of the event ${\cal E}(i,i')$ -- there is a path
connecting the two full slices.  The sites
$v_1,w_1,v_2,w_2,\ldots$ are shown as black dots in order along
the path from left to right.  In this example $v_3=w_3$, and
$v_4,v_6$ lie in the same slice, but in distinct components within
the slice.} \label{vw}
\end{figure}

The occurrence of ${\cal E}(i,i')$ implies the existence of
$v_1,\ldots ,w_s$ satisfying the above properties.  For $v,w\in
T_j$ define the event
\begin{align*}
\big\{v\conn{\leq n} w\big\}:=\Big\{\text{there exists a
connected,
$(d-1)$-internally-spanned set}\;\;& \\
\text{$S\subset T_j$ with $\diam\, S\leq n$ and $v,w\in S$}&\Big\}.
\end{align*}
(The definition of this event is delicate, and corrects a small
error in \cite{cerf-m}).  Then we have
\begin{align}
\nonumber \prob \big({\cal E}(i,i')\big)&\leq \sum_{s\geq
i'-i+1}\;\sum_{v_1,\ldots ,w_s}
\prob \Big(\{v_1\conn{\leq n}w_1\}\circ \cdots \circ\{v_s\conn{\leq n}w_s\}\Big) \\
\label{b-k} &\leq \sum_{s\geq i'-i+1}\;\sum_{v_1,\ldots ,w_s}
\prob \big(v_1\conn{\leq n}w_1\big) \cdots \prob\big(v_s\conn{\leq
n}w_s\big)
\end{align}
Here the second sum is over all possible choices of
$v_1,w_1,v_2,\ldots ,w_s$ satisfying properties (i)--(iii) above, the
symbol $\circ$ denotes disjoint occurrence (which holds because of
property (vi)), and the second inequality follows from the Van den
Berg-Kesten inequality (see e.g. \cite{g2}). In order to bound the
above, consider choosing $v_1,w_1,v_2,\ldots ,w_s$ in order. There are
$\#T_i=m^{d-1}$ possible choices for $v_1$.  Once $v_t$ is chosen, the
possible choices for $w_t$ lie in the cube
$$V(v_t):=v_t-(n+1,\dots ,n+1,0)+Q^{d-1}(2n+1)$$
centred at $v_t$.  Furthermore the event $\big\{v_t\conn{\leq
n}w_t\big\}$ is contained in the event $\big\{v_t\conn{\langle
X\cap V(v_t)\rangle} w_t\big\}$.  Once $w_t$ is chosen, there are
(at most) $2$ possible choices for $v_{t+1}$, corresponding to the
two neighbouring slices. Hence we obtain
\begin{align*}
\sum_{v_1,\ldots ,w_s}
\prob (v_1\conn{\leq n}w_1) \cdots \prob(v_s\conn{\leq n}w_s) \\
\leq m^{d-1} \bigg(\sum_{w\in V(0)}
 f_{2n+1}^{d-1}(0,w)\bigg)^s 2^{s-1}
\end{align*}
Substituting into (\ref{b-k}) we obtain
\begin{equation}
\label{bound-eii} \prob \big({\cal E}(i,i')\big)\leq 1\wedge m^{d-1}
H(i'-i),
\end{equation}
where $H(r)$ is as in the statement of Lemma \ref{constr}.

\sloppy We use an almost identical argument to bound the
probability of \linebreak ${\cal G}(x,i,i')$.  In this case the
path starts at $x$, so there is no need for the factor $m^{d-1}$.
\fussy We obtain
\begin{equation}
\label{bound-gii} \sum_a \prob \Big({\cal
G}(x,i,i')\cap\{A=a\}\Big)\leq H(i'-i).
\end{equation}

Finally, substituting
(\ref{prob-full}),(\ref{bound-eii}),(\ref{bound-gii}) into
(\ref{bound-1}),(\ref{bound-2}), and substituting these into
(\ref{case-1}),(\ref{case-2}) we obtain the conclusion of Lemma
\ref{constr}.
\end{pf}

In the following proofs we use $C_1,C_2,\ldots$ to denote
constants in $(0,\infty)$ which may depend on $d$ and $\epsilon$,
but not on $p$.

\begin{pfof}{Theorem \ref{sub-crit} (case $d=3$)}
Let $d=3$ and fix $\epsilon>0$.  Since $\epsilon$ is arbitrary we
can take for convenience
$$m=\Big\lfloor \exp{\frac{\lambda-2\epsilon}{p}}\Big\rfloor.$$
We shall bound $f^3_m(x,y)$ using Lemma \ref{constr}; for this we
need to choose $n$ and find upper bounds on $F^2_{m,n}$ and
$\chi_n^2$.

We first consider $F^2_{m,n}$.  If $\langle X \cap
Q^2(m)\rangle_2$ has a component of diameter $\geq n$ then by
Lemma \ref{aiz-leb}, $Q^2(m)$ contains a connected, 2-i.s.\ set
$T$ with diameter in $[n/2,n]$.  By Theorems 4 and 2(ii) of
\cite{h-boot} we may find $B=B(\epsilon)\in(0,\infty)$ such that
$I^2(B/p)\leq 1/\exp \frac{2\lambda-\epsilon}{p}$ \pss\ (note that
the factor of 2 in the exponent is important).  Indeed, by
equation (11) in \cite{h-boot} and the proof of Theorem 1(ii) in
\cite{h-boot}, we can choose $B$ sufficiently large that for any
connected $T\subset{\mathbb Z}^2$ with $\diam\,T\in\big[\lfloor
B/(2p)\rfloor,\lfloor B/p\rfloor\big]$ we have
$$\prob (T \text{ is 2-i.s.})\leq 1/\exp \frac{2\lambda
-2\epsilon}{p}
$$
\pss.  Therefore let
$$n=\Big\lfloor \frac{B}{p}\Big\rfloor.$$
Then by the above remarks we have
\begin{align}
\nonumber
F^2_{m,n}&\leq m^2 n^2 / \exp \frac{2\lambda-2\epsilon}{p} \\
\nonumber &\leq \bigg(\frac B p\bigg)^2
\exp\frac{2\lambda-4\epsilon-2\lambda+2\epsilon}{p} \\
\label{f-2} &\leq \exp\; -\frac\epsilon p
\end{align}
\pss.

We now turn to $\chi_n^2$, which is the expected volume of the
cluster at $(n+1,n+1)$ in $\langle X\cap Q^2(2n+1)\rangle$.  This
cluster is a rectangle, $R$ say.  $R$ is 2-i.s., and $\#R\leq
(\diam\, R)^2$. We bound $\chi_n^2=\expe_p \#R$ by considering two
cases. If $1\leq \diam\, R\leq 10$, then some site within distance
10 of $z$ must be occupied. If $\diam\, R>10$ then by Lemma
\ref{aiz-leb}, $Q^2(2n+1)$ must contain some 2-i.s.\ rectangle $S$
with diameter in $[5,10]$; and a 2-i.s.\ rectangle has at least
one occupied site in each row and each column. Hence we have
\begin{align}
\nonumber
\chi_n^2& \leq  10^2 (21^2 p) + (2n+1)^2 [(2n+1)^2 10^2] (10 p)^5 \\
\label{chi-2}
 &\leq C_1 p \leq \surd p,
\end{align}
\pss. (Here $(2n+1)^2 10^2$ is a bound on the number of possible
choices for the rectangle $S$, $(10p)^5$ is a bound on the
probability $S$ is internally spanned, and we have used the
definition of $n$).

Now we use (\ref{f-2}),(\ref{chi-2}) to bound the terms in Lemma
\ref{constr}.  We have \pss
$$H(r)\leq \frac 12\sum_{s=r-1}^\infty (2\surd p)^s \leq \frac{(2\surd
p)^{r-1}}{2(1-2\surd p)}\leq  p^{C_2(r-1)}.$$  Writing $C=C_2$ we
now bound the following expression from Lemma \ref{constr} by
considering two possible cases for the value of $r$:
\begin{eqnarray}
\nonumber
\lefteqn{ F_{m,n}^{2}[1\wedge m^2 H(r)]\leq  e^{- \epsilon /
p}\;\Big[1\wedge e^{2\lambda/p}\; p^{C(r-1)}\Big]} \\
\nonumber
&\leq &
\begin{cases}
 e^{-\epsilon/p}\leq p\,e^{-\epsilon/(2p)}
 \leq p \;p^{\textstyle[\frac{\epsilon C}{8\lambda}(r-1)]}
&\quad\text{if }p^{C(r-1)}\geq e^{-4\lambda/p}\\[2mm]
p  \Big[ e^{2\lambda/p}\; p^{C(r-1)/2}\;p^{C(r-1)/2}\Big]
 \leq p\; p^{C(r-1)/2}
 &\quad\text{if }p^{C(r-1)}< e^{-4\lambda/p}
\end{cases} \\
\label{cases-bound}
&\leq & p^{C_3 r}
\end{eqnarray}
\pss.
Looking again at Lemma \ref{constr} we therefore have
\begin{equation}
\label{funny-bound}
H(i_1+1)+\sum_{a=0}^\infty F_{m,n}^{2}\Big[1\wedge m^{2}
 H(i_1+a)\Big]
 \leq  p^{C_2 i_1} + C_4 p^{C_3 i_1} \leq p^{C_5 i_1}
\end{equation}
 \pss, where we can drop the initial multiplicative constant because
 $i_1\geq 1$.  Finally, substituting
 (\ref{cases-bound}),(\ref{funny-bound}) into Lemma \ref{constr} we
 obtain \pss
\begin{align*}
f_m^3(x,y)&\leq \sum_{k=0}^\ell \;\;
\sum_{0<i_1<\cdots<i_k<i_{k+1}=\ell+1}
 p^{C_5 i_1} p^{C_3 (i_2-i_1)} \cdots p^{C_3 (i_{k+1}-i_k)} \\
 &\leq 2^\ell p^{C_5(\ell+1)} \leq p^{\gamma'(\ell+1)}
\end{align*}
for some $\gamma'=\gamma'(d,\epsilon)>0$, as required.
\end{pfof}

\begin{pfof}{Theorem \ref{sub-crit} (case $d\geq 4$)}
The proof is by induction on dimension.  Fix $d\geq 4$ and
$\epsilon>0$ and suppose the case $d-1$ is proved.
Let
$$m=\Big\lfloor\exp^{d-2}{\frac{\lambda-\epsilon}{p}}\Big\rfloor
\quad\text{and}\quad
n=\frac 12 \Big\lfloor\exp^{d-3}{\frac{\lambda-\epsilon}{p}}\Big\rfloor-1,
$$
so that the inductive hypothesis gives for $x,y\in Q^{d-1}(n)$ that
$f_n^{d-1}(x,y)\leq f_{2n+1}^{d-1}(x,y)\leq p^{\gamma(\|x-y\|_\infty+1)}$.

We shall apply Lemma \ref{constr}.
If $\langle X \cap Q^{d-1}(m)\rangle_{d-1}$ has a component of
diameter $\geq n$ then by Lemma \ref{aiz-leb}, $Q^{d-1}(m)$
contains a connected, $(d-1)$-i.s.\ subset $T$ with diameter in
$[n/2,n]$.  Hence by the inductive hypothesis together with the
argument used to obtain (\ref{sub-set}) in the proof of Theorem
\ref{main}(ii), we have
$$F^{d-1}_{m,n} m^{d-1}\leq m^{2(d-1)} n^{2(d-1)} p^{\gamma (n/2+1)}\leq
1/\exp^{d-2}\frac{\lambda-\epsilon}{p}\leq p$$ \pss. (To check the
second inequality, take the logarithm).

Using the inductive hypothesis again we have \pss
$$\chi^{d-1}_n\leq \sum_{r=0}^n (2r+1)^{d-1} p^{\gamma (r+1)}
\leq p^\gamma \bigg(1+C_6\int_1^\infty r^{d-1}p^{\gamma r}
\,dr\bigg)\leq p^{C_7},$$
since \pss, $r^{d-1}p^{\gamma r}$ is decreasing in
 $r\geq 1$.  So \pss\ we have
$$H(r)\leq p^{C_8(r-1)}
\quad\text{and}\quad
 F^{d-1}_{m,n}\big[1\wedge m^{d-1}H(r)\big]\leq p^{C_9 r}.$$

Hence, as in the case $d=3$, substituting into Lemma \ref{constr}
gives \pss
$$f_m^d(x,y)\leq 2^\ell p^{C_9(\ell+1)}\leq p^{\gamma'(\ell+1)},$$
as required.
\end{pfof}

\section*{Open Problems}
\begin{mylist}
\item Prove the analogue of Theorem \ref{main} for the standard
bootstrap percolation model in $3$ or more dimensions.  What is
the value of the threshold $\lambda$ in this case?

\item Currently all proofs of the existence of a sharp threshold
(in the sense of Theorem \ref{main} as opposed to
\cite{balogh-bollobas-sharp}) for bootstrap percolation models
involve calculating its value. (See \cite{h-boot},
\cite{h-liggett-romik} and the present work). Is there a simpler
method of proving existence without determining the value?

\item What is the ``second order" asymptotic behaviour of the
model? Specifically, for example, if $p=p_{1/2}(L)$ is such that
$I^d(L,p)=1/2$, what is the asymptotic growth rate of $p\log^{d-1}
L - \lambda$ as $L\to\infty$?
\end{mylist}

\section*{Acknowledgement}
I thank Raphael Cerf for stimulating conversations.

\bibliography{mybib}

\end{document}